\documentclass[11pt]{article}

\usepackage{amsfonts}
\usepackage{amssymb,amsmath} 

\newcommand{\nicefrac}[2]
{\leavevmode \kern.1em\raise.5ex\hbox{\the\scriptfont0 #1}
             \kern-.1em/\kern-.15em\lower.25ex
             \hbox{\the\scriptfont0 #2}}

\newtheorem{Theo}{Theorem}{\alph{enumi}}
\newenvironment{theorem}{\begin{Theo}\hspace{-0.2cm}: }{\end{Theo}}

\newtheorem{Pro}{Proposition}{\alph{enumi}}

\newtheorem{Co}{Corollary}{\alph{enumi}}
\newenvironment{corollary}{\begin{Co}\hspace{-0.2cm}: }{\end{Co}}

\newtheorem{As}{Assumption}{\alph{enumi}}

\newtheorem{Le}{Lemma}{\alph{enumi}}
\newenvironment{lemma}{\begin{Le}\hspace{-0.2cm}: }{\end{Le}}

\newtheorem{Fo}{Folgerung}{\alph{enumi}}

\newtheorem{De}{Definition}{\alph{enumi}}

\newtheorem{Be}{Bemerkung}{\alph{enumi}}

\newtheorem{Ex}{Example}{\alph{enumi}}
\newenvironment{example}{\begin{Ex}\hspace{-0.2cm}:\rm }{\end{Ex}}

\newcommand{\ol}{\overline}
\newcommand{\bee}{\begin{equation}}
\newcommand{\ee}{\end{equation}}

\newcommand{\bea}{\begin{eqnarray}}
\newcommand{\ea}{\end{eqnarray}}

\newcommand{\ve}{\varepsilon}

\begin{document} 
\begin{center}
{\Large{\sc A uniqueness and periodicity result for solutions of elliptic 
equations in unbounded domains}}\\[2cm]
{\large{\sc Matthias Bergner, Jens Dittrich}}\\[1cm]
{\small\bf Abstract}\\[0.4cm]
\begin{minipage}[c][2cm][l]{12cm}
{\small We proof a uniqueness and periodicity theorem for bounded solutions of uniformly elliptic equations
in certain unbounded domains. } 
\end{minipage}
\end{center}
\subsection{Introduction}
In this note we study solutions $u\in C^2(\Omega,\mathbb R)\cap C^0(\ol\Omega,\mathbb R)$
of the Dirichlet problem
\bee\label{l1}
 a^{ij}(x) \partial_{ij} u+ b^i(x) \partial_i u+c(x) u=f
\quad \mbox{in}\; \Omega \quad, \quad u=g \quad \mbox{on}\; \partial\Omega \; . 
\ee
(using the sum convention) assuming the differential equation to be elliptic, i.e.
at each point $x\in\ol\Omega$ the matrix $a_{ij}(x)$ is symmetric and positive definite. 
In addition, we require the sign condition $c(x)\leq 0$.\\ \\
%A linear elliptic problem of type (\ref{l1}) arises when studying the difference function
%$u(x)=u_1(x)-u_2(x)$ of two solutions $u_1,u_2$ of some quasilinear, elliptic problem
%\bee\label{l6} 
%A^{ij}(x,u_k,D u_k)\partial_{ij} u_k+ B(x,u_k, Du_k)=0 \quad \mbox{in}\; \Omega
%\quad , \quad u_k(x)=g(x) \quad \mbox{on}\; \partial\Omega \; . 
%\ee
%If one is interested in uniqueness of solutions to the quasilinear equation (\ref{l6})
%one simply has to show that the corresponding linear equation (\ref{l1}) only admitts the
%trivial solution $u\equiv 0$. 
If the domain $\Omega$ is bounded,
the well known classical maximum principle (see \cite[Theorem 3.3]{gilbarg})
asserts that (\ref{l1}) admitts at most one solution.
In contrast, such a result
does in general not hold for unbounded domains $\Omega$ and examples are given below.
First, let us make the following assumptions
on the coefficients: Let $a^{ij},b^i,c\in C^0(\ol\Omega,\mathbb R)$ and satisfy
\bee\label{l11}
||a^{ij}||_{C^0(\Omega)}+||b^i||_{C^0(\Omega)}
+||c||_{C^0(\Omega)}\leq H \quad \mbox{for}\; i,j=1,\dots,n \quad \mbox{and} 
\quad c(x)\leq 0 \quad \mbox{in}\; \Omega 
\ee
with some constant $H$. Additionally,
we have to require a uniform ellipticity condition
\bee\label{l12}
\frac{1}{\Lambda}|\xi|^2\leq a^{ij}(x)\xi_i\xi_j\leq \Lambda |\xi|^2 
\quad \mbox{for all}\; x\in\Omega\; , \; \xi\in\mathbb R^n 
\ee
with constant $\Lambda<\infty$.
Our first result is the following
\begin{theorem}\label{t1}
Additionally to (\ref{l11}) and (\ref{l12}), assume the following
\begin{itemize}
\item[a)] The unbounded domain $\Omega\subset\mathbb R^n$ has bounded thickness, 
i.e. $\sup_{x\in\Omega}\mbox{dist}(x,\partial\Omega)<+\infty$.
\item[b)] Let $u\in C^2(\Omega,\mathbb R)\cap C^0(\ol\Omega,\mathbb R)$
be a bounded solution of (\ref{l1}) for the right side $f\equiv 0$ and boundary values $g\equiv 0$.
\item[c)] Let $u$ satisfy the following uniform boundary condition:
For any sequence $x_k\in\Omega$ with $\mbox{dist}(x_k,\partial\Omega)\to 0$
for $k\to\infty$ it follows that $u(x_k)\to 0$ as $k\to\infty$.
\end{itemize}
Then we must have $u\equiv 0$ in $\Omega$.
\end{theorem}
As the proof of Theorem \ref{t1} reveals, this result remains true
for weak solutions $u$ of regularity class 
$W^{2,n}_{loc}(\Omega,\mathbb R)\cap C^0(\ol\Omega,\mathbb R)$. \\ \\
Let us now demonstrate the necessity of the assumption a), b) and c)
by considering the following examples.
\begin{example}
For some $k\in\mathbb N$ take the domain
$\Omega=\{re^{i\varphi}\in\mathbb C \, |\, 0<\varphi<\frac{\pi}{k}\}$
and the harmonic function $u(x,y):=\mbox{Re}\{(x+iy)^k\}$
with $u=0$ on $\partial\Omega$. This very simple example already shows
that certain assumption on the solution $u$ or the domain $\Omega$
are needed for a uniqueness theorem to hold. Note that for this example
all of the assumptions a), b) and c) of Theorem \ref{t1} are not satisfied.
\end{example}
\begin{example}
As domain we take $\Omega=\{(x,y)\in\mathbb R^2\, | \, 0<y<\pi\}$
and consider the unbounded, harmonic function $u(x,y)=e^x \sin y$ with $u=0$ on $\partial\Omega$.
Here, assumption a) of Theorem \ref{t1} is satisfied while assumptions b) and c) are not.
\end{example}
\begin{example}
Now consider the domain $\Omega=\{x\in\mathbb R^n\, : \, |x|>1\}$, $n\geq 3$
and the bounded, harmonic function $u(x)=1-|x|^{2-n}$ with $u=0$ on $\partial\Omega$.
Here, assumptions b) and c) of Theorem \ref{t1} are satisfied while assumption a) is not.
\end{example}
Let us make a remark on assumption c): 
If the domain $\Omega$ had a compact boundary $\partial\Omega$, 
then assumption c) would directly follow from 
$u\in C^0(\ol\Omega,\mathbb R)$ together with $u=0$ on $\partial\Omega$.
However, note that an unbounded domain cannot both have a compact
boundary and at the same time satisfy assumption a).
Assumption c) will hold provided that the solution $u$ is uniformly continuous in $\Omega$.
By suitably restricting the domain $\Omega$, we can actually show a uniform
continuity of the solution.
\begin{theorem}\label{t2}
Additionally to (\ref{l11}) and (\ref{l12}), assume the following
\begin{itemize}
\item[a)] The unbounded domain $\Omega\subset\mathbb R^n$ has bounded thickness, 
i.e. $\sup_{x\in\Omega}\mbox{dist}(x,\partial\Omega)<+\infty$.
\item[b)] Let $u_k\in C^2(\Omega,\mathbb R)\cap C^0(\ol\Omega,\mathbb R)$
be two solutions of (\ref{l1}) for some right side $f$ and some boundary values $g$.
Assume difference $|u_1(x)-u_2(x)|$ is uniformly bounded in $\Omega$.
\item[c)] The domain $\Omega$ satisfies a uniform exterior sphere condition.
\end{itemize}
Then it follows $u_1\equiv u_2$ in $\Omega$. 
\end{theorem}
By uniform exterior sphere condition we mean the following: There exists some
$r>0$ such that for each $x_0\in\partial\Omega$ there exists
some $x_1\in\mathbb R^n$ such that $\ol{B_r(x_1)}\cap\ol\Omega=\{x_0\}$. \\ \\
Finally, we want to point out that uniqueness results for partial
differential equations also imply symmetry properties of the solutions.
To illustrate this by an example, we have the following result.
\begin{corollary}\label{c1}
Assume that (\ref{l11}) and (\ref{l12}) hold.
Moreover, assume that $\Omega$ satisfies a uniform exterior 
sphere condition and can be decomposed into
$\Omega=\mathbb R\times\Omega'$ for some bounded domain
$\Omega'\subset\mathbb R^{n-1}$.
We require the coefficients $a^{ij}$, 
$b^i$, $c$, the right side $f$ and the boundary values $g$ to be periodic
w.r.t. the $x_1$-variable with one and the same period length $L>0$. \\
Then any bounded solution $u\in C^2(\Omega,\mathbb R)\cap C^0(\ol\Omega,\mathbb R)$
of (\ref{l1}) is periodic w.r.t. the $x_1$-variable.
\end{corollary}
Note that in Example 2 we found a solution not being periodic.
There, all of the assumptions of Corollary \ref{c1} are satisfied
except for the boundedness of the solution. Hence, also for
Corollary \ref{c1} it is crucial only to consider bounded solutions. 

%%%%%%%%%%%%%%%%%%%%%%%%%%%%%%%%%%%%%%%%%%%%%%%%%%%%%%%%%%%%%%%%%%%%%%%%%%%%%%
\subsection{The proof of Theorem \ref{t1}}
For the proof of Theorem \ref{t1}, we first need the following lemma,
which may be of independent interest. It is a generalisation of the strong maximum principle.
\begin{lemma}\label{lemma1}
Let $u_k\in C^2(\Omega,\mathbb R)$ be a sequence of solutions of 
$$ a^{ij}_k(x) \partial_{ij} u_k
+b^i_k(x)\partial_i u_k+c_k(x)u_k =0\quad \mbox{in} \; \Omega \; . $$
Let the coefficients $a^{ij}_k,b^i_k,c_k$ satisfy
(\ref{l11}) and (\ref{l12}) with constants $\Lambda,H$ independent of $k$ and $c_k(x)\leq 0$ in $\Omega$.
Assume that $u_k$ converge uniformly in $\Omega$ to some $u\in C^0(\Omega,\mathbb R)$.
For some $x_*\in\Omega$ and $M\in\mathbb R$ let
$$ u_k(x)\leq M \quad \mbox{in} \; \Omega \quad \mbox{and} \quad \lim\limits_{k\to\infty} u_k(x_*)=M \; . $$
Then it follows $u\equiv M$ in $\Omega$.
\end{lemma}
{\it Proof:} \\
Consider the set
$$ \Theta:=\{x\in \Omega\, | \, u(x)=M\} $$
which is not empty because of $x_*\in \Theta$. Now $\Theta$ is closed within $\Omega$ due to the continuity of $u$.
We now show that $\Theta$ is also open implying $\Theta=\Omega$ and proving the lemma.
For $x_0\in \Theta$ choose $r>0$ small enough such that $\ol{B_{2r}(x_0)}\subset \Omega$.
Now consider the function $v_k(x):=M-u_k(x)$ for $x\in\Omega$ with $v_k(x)\geq 0$ in $\Omega$.
Because of $c_k\leq 0$ this $v_k$ is a solution of the differential inequality
$$ a^{ij}_k(x) \partial_{ij} v_k +b^i_k(x)\partial_i v_k\leq 0 \quad \mbox{in}\; \Omega \; . $$
We now apply the Harnack type inequality \cite[Theorem 9.22]{gilbarg} on the
domain $B_{2r}(x_0)$: There exist constants $p>0$ and $C<\infty$
only depending on $r,H,\Lambda$ and $n$ such that
\bee\label{l20}
\Big\{\int\limits_{B_r(x_0)} v_k(x)^p dx\Big\}^{1/p}\leq 
C \inf\limits_{B_r(x_0)} v_k(x)\leq C v_k(x_0)=C\Big(M-u_k(x_0)\Big) \; .
\ee
Noting that $u_k(x_0)\to u(x_0)$ for $k\to\infty$ and $u(x_0)=M$ because of $x_0\in\Theta$,
passing to the limit in (\ref{l20}) then yields
$$ \Big\{\int\limits_{B_r(x_0)} (M-u(x))^p dx\Big\}^{1/p}= 0 \; . $$ 
Together with $u(x)\leq M$ in $\Omega$ this implies 
$u(x)=M$ in $B_r(x_0)$ proving that $\Theta$ is open. \hfill $\Box$ \\ \\
{\it Remarks:} 
\begin{itemize}
\item[1.)] The lemma remains true for weak solutions of
regularity $W^{2,n}(\Omega,\mathbb R)\cap C^0(\Omega,\mathbb R)$.
\item[2.)]
The proof of this lemma is similar to the proof
of the strong maximum principle for weak solutions (see \cite[Theorem 8.19]{gilbarg}).
In case of $u_k(x)=u(x)$ for all $k$ the statement of the lemma 
reduces to the classical strong maximum principle. 
\end{itemize}
{\it Proof of Theorem \ref{t1}:} \\
Given a solution $u$ of (\ref{l1}) for $f\equiv 0$ and $g\equiv 0$,
we will show $u\equiv 0$ in $\Omega$ as follows:
Assume to the contrary that $u(x_0)\neq 0$ for some $x_0\in\Omega$, 
say $u(x_0)>0$. Defining $M:=\sup_\Omega u(x)>0$
we have $M<+\infty$ by the boundedness assumption b) of Theorem \ref{t1}.
We can now find a sequence $x_k\in\Omega$ such that $u(x_k)\to M$ for $n\to\infty$.
%Two cases are possible: \\ \\
%1st case: The sequence $x_k$ contains a bounded subsequence.
%After extracting a convergent subsequence we may then assume that $x_k\to x_*$ for $n\to\infty$ with some $x_*\in \ol\Omega$.
%By the continuity of $u$ we conclude $u(x_k)\to u(x_*)=M$ for $k\to\infty$.
%Together with $u=0$ on $\partial\Omega$ this implies $x_*\in\Omega$.
%Thus the function $u$ achieves its maximum at the point $x_*\in\Omega$
%and by the strong maximum principle (see \cite[Theorem 3.5]{gilbarg})
%we have $u(x)\equiv u(x_*)=M$ in $\Omega$, contradicting $u=0$ on $\partial\Omega$. \\ \\
Now, for each $k\in\mathbb N$ let us define
$$r_k:=\mbox{dist}(x_k,\partial\Omega)  \; .$$
We claim that there exist constants $\ve>0$ and $R<\infty$ such that
\bee\label{l3}
\ve<r_k<R \quad \mbox{for all}\; k\in\mathbb N \; . 
\ee
In fact, the right inequality follows directly from assumption a) of Theorem \ref{t1}
if we define $R=\sup_{x\in\Omega} \mbox{dist}(x,\partial\Omega)$.
The left inequality follows from assumption c) together with $u(x_k)\to M>0$. 
On the ball $B:=\{x\in\mathbb R^n\, : \, |x|<1\}$,
let us now consider the shifted and rescaled functions
$$ v_k:\ol B\to\mathbb R \quad , \quad v_k(x):=u(x_k+r_k x) \; . $$
By the definiton of $r_k$, for each $k\in\mathbb N$ we can find some $y_k\in\partial B$
with $x_k+r_k y_k\in\partial\Omega$ implying $v_k(y_k)=u(x_k+r_k y_k)=0$.
Since $u$ is solution of (\ref{l1}), $v_k$ will then be solution of
$$ a^{ij}_k(x) \partial_{ij} v_k+b^i_k(x) \partial_i v_k+c_k(x) v_k=0
\quad \mbox{in}\; B \; \mbox{for}\; k\in\mathbb N $$
with coefficients $a^{ij}_k(x):=r_k^{-2}a^{ij}(x_k+r_k x)$ and $b^i_k$, $c_k$ defined similarly.
By (\ref{l3}) together with the assumptions on $a^{ij},b^i,c$
there is a uniform $C^0$-bound
$$ \sup\limits_{k\in\mathbb N}
\Big(||a^{ij}_k||_{C^0(B)}+||b^i_k||_{C^0(B)}+||c_k||_{C^0(B)}\Big)<+\infty
\quad \mbox{for all}\; i,j=1,\dots,n \; . $$
Using the interior H\"older estimate \cite[Theorem 9.26]{gilbarg} for weak solutions we get
$$ \sup\limits_{k\in\mathbb N} ||v_k||_{C^\alpha(B_s)}<+\infty \quad \mbox{for all} \; 0<s<1 $$
with some H\"older exponent $\alpha=\alpha(s)\in (0,1)$ independent of $k$.
After extracting some subsequence we obtain the uniform convergence
\bee\label{l4}
 v_k \; \to \; v \quad \mbox{in} \; C^0(B_s,\mathbb R) \; \mbox{for} \; k\to\infty 
\ee
for each $s<1$ with some limit function $v\in C^0(B,\mathbb R)$ satisfying
$$ v(x)\leq M \quad \mbox{in} \; B \quad \mbox{and}\quad v(0)=M \; . $$
By Lemma \ref{lemma1} (applied to $\Omega=B_s$) we have 
$v(x)=M$ in $B_s$ for each $s<1$ and hence $v(x)=M$ in $B$. \\ \\
On the other hand, from $v_k(y_k)=0$ together with $v_k(0)\to M$ we conclude that, for
sufficiently large $k$, there exists some
$z_k=t_k y_k\in B$ with $t_k\in (0,1)$ such that
$v_k(z_k)=M/2=u(x_k+r_k z_k)$. We may assume that
$t_k\to t_*\in [0,1]$ and $z_k\to z_*\in\ol B$ as $k\to\infty$.
We now claim that $t_*<1$. Otherwise we would have $t_k\to 1$ for $k\to\infty$.
However, we would then have 
$$ \mbox{dist}(x_k+r_k z_k,\partial\Omega)
\leq |x_k+r_k z_k-(x_k+r_k y_k)|=r_k |y_k|(1-t_k)\leq R (1-t_k)\to 0 \quad \mbox{for}\; k\to\infty$$
contradicting assumption c) together with $u(x_k+r_k z_k)=M/2$, proving the claim.
Using the uniform convergence (\ref{l4}) in the ball $B_{t_*}$
together with $M/2=v_k(z_k)$ 
we obtain $v(z_*)=M/2$, contradicting $v(x)\equiv M$ in $B$. \hfill $\Box$
\subsection{The proof of Theorem \ref{t2} and Corollary \ref{c1}}
We start with \\ \\
{\it Proof of Theorem 2:} \\
Consider two bounded solutions $u_1,u_2$ of (\ref{l1}). Then the difference
function $u(x):=u_1(x)-u_2(x)$ will be solution of (\ref{l1}) for the right
side $f\equiv 0$ and boundary values $u\equiv 0$ on $\partial\Omega$.
By assumption b) of Theorem \ref{t2}, $u$ is bounded in $\Omega$,
hence $|u(x)|\leq M$ for some $M>0$.
We want to apply Theorem \ref{t1} to $u$, but we first have to check
wether the uniform boundary condition, assumption c) of Theorem 1, is satisfied by $u$. 
As described in Remark 3 of \cite[Chapter 6.3]{gilbarg}
we can construct a uniform barrier at each boundary point, using the uniform
exterior sphere condition.
Let $R>0$ be the radius of the uniform exterior sphere condition
and $x_0\in\partial\Omega$. Then there exists some $y\in\mathbb R^n$
with $\ol{B_R(y)}\cap \ol\Omega=\{x_0\}$.
Consider the function
$$ w(x):=R^{-\sigma}-|x-y|^{-\sigma} \quad \mbox{for}\; x\in\mathbb R^n \; , \; \sigma>0 $$
satisfying $w(x_0)=0$ and $w(x)>0$ in $\Omega$.
Setting $r:=|x-y|$ and using $c\leq 0$ in $\Omega$ we estimate
\bea
&& a^{ij}(x)\partial_{ij} w+b^i(x)\partial_i w+c(x) w \nonumber \\
&\leq&a^{ij}[-\sigma(\sigma+2) r^{-\sigma-4} (x_i-y_i) (x_j-y_i)+\delta_{ij}\sigma r^{-\sigma-2}]
+b^i \sigma r^{-\sigma-2} (x_i-y_i) \nonumber \\
&=&\sigma r^{-\sigma-4}[-(\sigma+2)a^{ij} (x_i-y_i) (x_j-y_j)+r^2(a^{ij} \delta_{ij}+b^i (x_i-y_i))] \nonumber \\
&\leq&\sigma r^{-\sigma-2}[-(\sigma+2)\Lambda^{-1}+n H+n(R+1) H] \nonumber
\ea
for all $x\in\tilde\Omega:=\{x\in\Omega : |x-y|<R+1\}$.
By choosing $\sigma=\sigma(\Lambda,H,n)>0$ sufficiently large, we obtain
$$ a^{ij}(x)\partial_{ij} w+b^i(x)\partial_i w+c(x) w\leq 0 \quad \mbox{in}\; \tilde\Omega \; . $$
We now define
$$ \tau:=\frac{M}{R^{-\sigma}-(R+1)^{-\sigma}}>0 $$
and note that $-\tau w(x)\leq u(x)\leq \tau w(x)$ on $\partial\tilde\Omega$.
From the maximum principle we conclude that $-\sigma w(x)\leq u(x)\leq \sigma w(x)$ in $\tilde\Omega$.
Using $|x_0-y|=R$ this yields
\bea
 |u(x)|&\leq& \tau |w(x)|=\tau \Big( R^{-\sigma}-|x-y|^{-\sigma}\Big) \nonumber \\
&\leq& \tau \Big(R^{-\sigma}-(|x-x_0|+R)^{-\sigma}\Big)
\quad \mbox{for all}\; x\in\Omega\; , \; x_0\in\partial\Omega \; \mbox{with}\; |x-x_0|<1 \; . \nonumber
\ea
In particular, for $|x_0-x|=\mbox{dist}(x,\partial\Omega)$ we obtain
$$ |u(x)|\leq \tau \Big(R^{-\sigma}-(\mbox{dist}(x,\partial\Omega)+R)^{-\sigma}\Big)
\quad \mbox{for all}\; x\in\Omega \; \mbox{with} \; \mbox{dist}(x,\partial\Omega)<1 \; . $$
As the constants $R,\sigma$ and $\tau$ are independent of the choosen
boundary point $x_0\in\partial\Omega$, we see that assumption c) of Theorem \ref{t1} 
is satisfied by $u$. \hfill $\Box$ \\ \\
We finally give the\\ \\
{\it Proof of Corollary \ref{c1}:} \\
Let $\Omega=\mathbb R\times\Omega'$ for some bounded domain $\Omega'\subset\mathbb R^{n-1}$.
Note that such a domain $\Omega$ satisfies the uniform thickness condition
$\sup_\Omega \mbox{dist}(x,\partial\Omega)\leq d$ with $d:=\mbox{diam}(\Omega')$.
Let $u\in C^2(\Omega,\mathbb R)\cap C^0(\ol\Omega,\mathbb R)$
be a bounded solution of (\ref{l1}). For some $k\in\mathbb Z$ let us define a 
translation of $u$ by
$$ \tilde u(x)\in C^2(\Omega,\mathbb R)\cap C^0(\ol\Omega,\mathbb R) \quad , \quad \tilde 
u(x_1,\dots,x_n):=u(x_1+k L,x_2,\dots,x_n) \quad \mbox{for} \; x\in\ol\Omega \; . $$
Note that $\tilde u$ is bounded just as $u$ is.
By the periodicity assumptions on the data $a^{ij}$, $b^i$, $c$, $f$ and $g$
this $\tilde u$ will be solution of the same problem (\ref{l1}) as $u$.
By Theorem \ref{t2} we obtain $\tilde u(x)=u(x)$ in $\ol\Omega$
proving the periodicity of $u$. \hfill $\Box$

\vspace*{1cm}
Matthias Bergner, Jens Dittrich\\
Universit\"at Ulm\\
Fakult\"at f\"ur Mathematik und Wirtschaftswissenschaften\\
Institut f\"ur Analysis\\
Helmholtzstr. 18\\
D-89069 Ulm \\
Germany\\[0.3cm]
e-mail: matthias.bergner@uni-ulm.de, jens.dittrich@uni-ulm.de\\

\begin{thebibliography}{10}
\bibitem{gilbarg}
D.Gilbarg, N.S.Trudinger: {\it Elliptic Partial Differential Equations of Second Order}. 
Springer, Berlin Heidelberg New York, 1983.
\bibitem{saubuch1}
F.Sauvigny: {\it Partielle Differentialgleichungen der Geometrie und der Physik, Teil 1 und 2}. 
Springer Berlin Heidelberg, 2004, 2005.
\bibitem{saubuch3}
F.Sauvigny: {\it Partial Differential Equations, Vol. 1 and 2}. Springer Universitext, 2006.
\end{thebibliography}
\end{document}